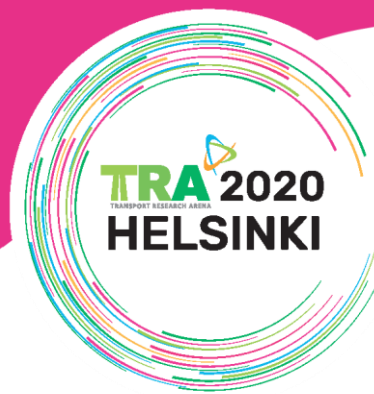

*Proceedings of 8th Transport Research Arena TRA 2020, April 27-30, 2020, Helsinki, Finland*

# aColor: Mechatronics, Machine Learning, and Communications in an Unmanned Surface Vehicle

Jose Villa[a]*, Jussi Taipalmaa[b], Mikhail Gerasimenko[c], Alexander Pyattaev[c], Mikko Ukonaho[d], Honglei Zhang[b], Jenni Raitoharju[b], Nikolaos Passalis[b], Antti Perttula[d], Jussi Aaltonen[a], Sergey Andreev[c], Markus Aho[d], Sauli Virta[e], Moncef Gabbouj[b], Mikko Valkama[c], Kari T. Koskinen[a]

[a]*Unit of Automation Technology and Mechanical Engineering, Tampere University (TAU), Korkeakoulunkatu 1, Tampere, 33720, Finland*
[b]*Unit of Computing Sciences, Tampere University (TAU), Korkeakoulunkatu 1, Tampere, 33720, Finland*
[c]*Unit of Electrical Engineering, Tampere University (TAU), Korkeakoulunkatu 1, Tampere, 33720, Finland*
[d]*School of Industrial Engineering, Tampere University of Applied Science (TAMK), Kuntokatu 3, Tampere, 33520, Finland*
[e]*R&D Department, Alamarin-Jet Oy, Tuomisentie 16, Härmä, 62300, Finland*

**Abstract**

The aim of this work is to offer an overview of the research questions, solutions, and challenges faced by the project aColor ("Autonomous and Collaborative Offshore Robotics"). This initiative incorporates three different research areas, namely, mechatronics, machine learning, and communications. It is implemented in an autonomous offshore multicomponent robotic system having an Unmanned Surface Vehicle (USV) as its main subsystem. Our results across the three areas of work are systematically outlined in this paper by demonstrating the advantages and capabilities of the proposed system for different Guidance, Navigation, and Control missions, as well as for the high-speed and long-range bidirectional connectivity purposes across all autonomous subsystems. Challenges for the future are also identified by this study, thus offering an outline for the next steps of the aColor project.

* Corresponding author. Tel.: +358 50 448 1926;
 *E-mail address:* jose.villa@tuni.fi

## 1. Introduction

Developing autonomous surface vehicles that can function without any human intervention can improve numerous offshore operations by making them more efficient and reliable, as well as safer. In this paper, we offer an overview of research questions, solutions, and challenges faced during the project aColor ("Autonomous and Collaborative Offshore Robotics") (TechFinland100, 2018). aColor aims to develop the basic methodology for autonomous offshore robotic systems, as well as to demonstrate them in challenging open environments. To this end, several aspects of the system architecture have been considered, ranging from structural and mechatronic design, sensing, communication, and task planning to perception, data collection, and augmented shared intelligence. While the project's main focus is set on Unmanned Surface Vehicles (USVs), the proposed autonomous robotic system also includes cooperative features between the USVs, an Unmanned Aerial Vehicle (UAV), and an Autonomous Underwater Vehicle (AUV). Hence, the project proposes an innovative approach for multipurpose offshore intervention tasks from communications and connectivity to obstacle detection and path planning. The connectivity between these modules involves many different research areas and it is critical for ensuring the desired safe and seamless operation of the USV. This becomes even more challenging in co-operative robotics, where multiple vehicles are collaborating to form a diverse robotic net that works as a unified system to accomplish certain tasks.

Guidance, Navigation, and Control (GNC) systems are among the most vital components of a USV and their modules are generally constituted by on-board computers and software, which combined are responsible for managing the entire USV system. Trajectory tracking and path following are among the fundamental research challenges in the control domain. Several studies include different guidance laws for path following and collision avoidance (Fossen, et al., 2003) (Moe & Pettersen, 2016). An appropriate GNC technique regarding its architecture is essential to efficiently model and design the complete motion control system.

Observation and autonomous detection of objects on the water surface is a key task to ensure the safe and efficient operation of an USV. Information about objects on the water surface can be collected using different types of equipment like RGB cameras, thermal cameras, marine radars, and LiDARs. Across the variety of sensing devices, RGB cameras have been extensively used in various autonomous driving systems due to their low cost and ability to provide high-resolution images. With the visual information acquired by cameras, objects on the water surface can be detected and recognized using various machine learning and computer vision techniques, allowing for making informed decisions promptly.

In order to provide uninterrupted data collection and the capability to control all subsystems of the robotic platform, it is important to enable reliable high-speed long-range communications among the core components of the system. While it is possible to implement some of those qualities using industrial or military-grade equipment, for the operation of a scalable unmanned marine ecosystem it is essential to employ expensive, commercially-available solutions. For this purpose, high-directional antenna systems can be used with several software and hardware modifications. Finally, some of the ecosystem elements, such as an UAV, can be introduced to further extend the coverage and capacity of the dynamic robotic network. Design and implementation of the software and hardware aspects of this network is one of the core challenges in the communications research area.

An automated landing and charging platform for an UAV was selected as an additional research goal for this study. With GPS, an UAV navigates in regular flight mode with the inaccuracy of position of around 1.5 meters. To achieve better accuracy, the UAV can use satellite navigation Real-Time Kinematic (RTK) mode (Eriksson, 2016), normal camera, or IR camera technology. For battery charging purposes, however, the position inaccuracy needs to be in the range of cm, for which an external electromechanical positioning solution has been developed.

The main contribution of this paper is to provide a systematic description of a complex autonomous USV platform, which was developed within the scope of the aColor project and comprises of several dissimilar research areas. The overall USV platform is introduced in Section 2.1-2.4. The proposed system is further studied from three different points of view: mechatronics (Section 2.6), machine learning (Section 2.7), and communications (Section 2.8). Research and key challenges are described for all these fields as well as in combining them into a unified system. Furthermore, aerial vehicles architecture and a landing platform for the UAV located at the USV are described to demonstrate how cooperative robotics can be employed to further improve the working capabilities of the proposed offshore system (Section 3).

## 2. USV main platform

*2.1. Description of the USV platform*

The vessel utilized in this research comprises of an aluminium hull with thrust vectoring twin water-jet configuration. It uses two marine diesel engines with 170 kW of rated power and a control system developed by Alamarin-Jet Oy. The research vessel additionally contains a built-in interface for remote and autonomous operations, which can be used for research purposes in GNC applications. Furthermore, the aColor USV has an excellent manoeuvring accuracy, as it can move in all the directions in the horizontal plane without bow thrusters.

*2.2. Instrumentation of the USV platform*

The USV employed in this study carries various instrumentation that can be used in a wide variety of different applications. These include control (system and mission), situational awareness (including object detection and object recognition), and high-speed wireless connectivity. The complete instrumentation placed in the aColor USV is detailed in Fig. 1(a).

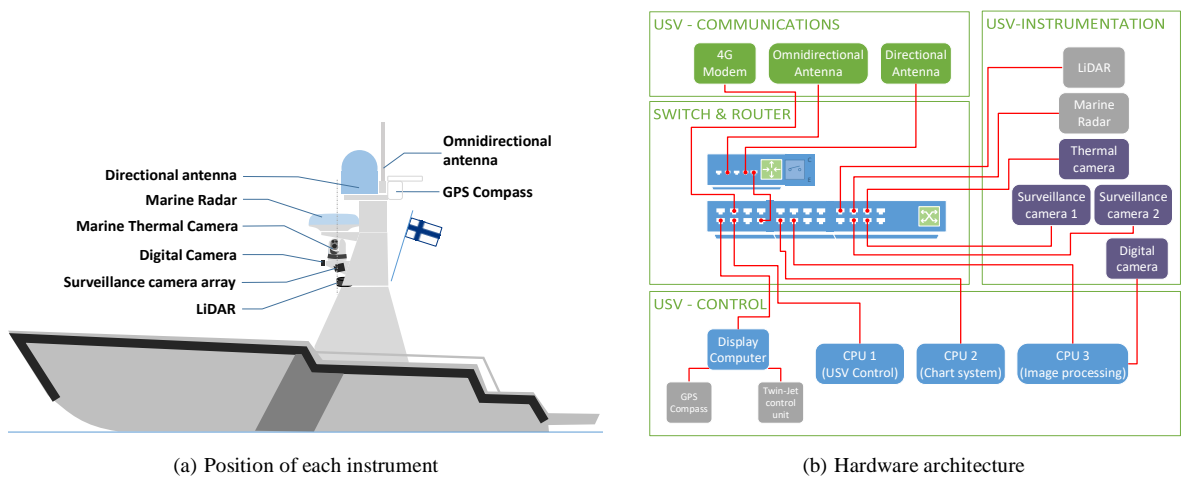

(a) Position of each instrument   (b) Hardware architecture

Fig. 1 aColor USV instrumentation

*2.3. Hardware schematic of the USV platform*

The complete hardware architecture utilized in the aColor USV is shown in Fig. 1(b), which also highlights the three main research areas: communications, situational awareness, and control. USV communications (green blocks) feature 4G along with omnidirectional and directional connections that can be used depending on the scenario. USV instrumentation can be generally grouped into two different categories of environment perception: active perception methods (grey blocks) and passive perception methods (purple blocks). LiDAR and marine radar are employed for obstacle avoidance, short and long range, respectively. Thermal camera (infrared vision), surveillance cameras (stereo vision), and digital camera (monocular vision with a wide field of view) are used in the USV for obstacle recognition capabilities. The aim of introducing them is to complement the information obtained from the active perception methods by identifying safe objects (e.g., a seagull) and by enabling active planning (e.g., calculating the expected path for another vessel and then planning the path of the USV accordingly). The USV control consists of different computers, which are responsible for sending and receiving commands related to GNC operations as well as for processing the data collected from various on-board sensors.

*2.4. Mast structure optimization*

The original design of the USV mast structure, developed without mechanical analyses, was found to be prone to vibrations. Operating the USV on wavy waters caused mechanical shocks to the mast structure and thus severe vibration, especially harmful for the directional antenna. Re-design of the mast was started by vibration measurements using vibration input module and accelerometers. Measured data was analysed to produce both magnitude and frequency profile of the vibration. Measured data was compared against simulation data from FEM-simulation. One of the eigenfrequencies in simulation results was very close to the dominating frequency in the measured data. The corresponding eigenmode showed that the mast was mainly moving in horizontal direction. Transient response of the mast was simulated using FEM-model and different design choices were analysed and compared against the original design. Finally, topology optimization with stiffness objective and volume constraint was applied to determine the optimum structure (solidThinking Inspire, 2018), (Turkkila & Ukonaho, 2016). As a result, the shape of the mast structure was fixed by using aluminium sheet metal to manufacture the parts as well as the placement of the radar. Fig. 2 presents the results of our modal analysis for both the original and the modified structure. Lower part of the mast was made stiffer by adding plexiglass windshield to the structure. Top part of the mast was completely re-designed, while stiffness in the horizontal direction was optimized against dynamic loadings. With re-designed parts, the lowest eigenfrequency of the mast is roughly 2.5 times higher and displacements under transient loading are two times lower as compared to the original design.

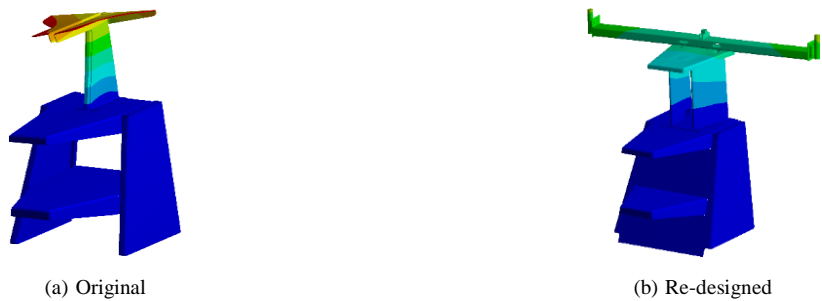

(a) Original  (b) Re-designed

Fig. 2 USV mast structures

*2.5. Collaboration & challenges in real demonstration*

A series of public demonstrations was conducted in the Tampere area (Finland) to confirm the needed capabilities of our developed autonomous system. These capabilities include the main research areas focused in this study: GNC algorithms related to the control architecture of the subsystem (mechatronics), machine learning algorithms for obstacle avoidance and object recognition, and long-range communications between the Ground Control Station (GCS) and the USV.

There are certain challenges that need to be considered in the study, such as the possibility of losing connectivity between the GCS and the USV, which can be critical in practical situations. Complete autonomous operation is the main challenge for both machine learning algorithms (data collection) and GNC testing algorithms. Furthermore, the position and orientation requirements for the instrumentation in the USV mast have been a crucial part in the development of every research topic, since the equipment for communication, collision avoidance, and obstacle detection requires a certain position to have its full working capability. Cameras, radar, and directional antenna have been the most important components of our USV system.

*2.6. USV GNC architecture*

*2.6.1. Research goals*

The key research objective in the mechatronics area is the design of a model-based architecture, which can be used for both simulation and implementation in different subsystems. Our architecture contains the main GNC methods, including control, navigation, communication, and data processing for any sensor featured by the network, while at same time providing the capability to easily add or remove any component of the system. Applying the same architecture to all of the offshore vehicles makes it possible to combine independent capabilities of the individual instrumentation into a more holistic result (e.g., a 3D point-cloud map for collision avoidance purposes).

*2.6.2. Methods and results*

The USV used in this study contains a Linux computer as a Master node for the system (CPU1), which is connected to the rest of the instrumentation by a switch, as it is shown in Fig. 1(b). This Linux computer has ROS installed to send and receive the necessary commands for the USV operation (Quigley, et al., 2019). The display computer connects the GPS compass and the twin-jet control unit with the rest of the system. This computer is used to receive the CAN messages from the GPS compass; it also sends the necessary commands to the twin-jet control unit, generating a specific thrust vector for each movement produced by the control unit. LiDAR and marine radar are included in the system for collision avoidance capabilities of the USV as well as to construct a 3D map of the environment. All this instrumentation installed in the USV employs ROS as a framework. Hence, this framework allows the use of the necessary tools and packages to access sensor data, process it, and generate an appropriate response for different actuators.

Fig. 3 displays the described schematic used for GNC methods, where the MATLAB/Simulink computer is only used for testing. A standalone ROS-node (The Mathworks, Inc., 2018) permits a rapid prototyping procedure while testing, and it can be run directly in the Linux computer which operates as a ROS master.

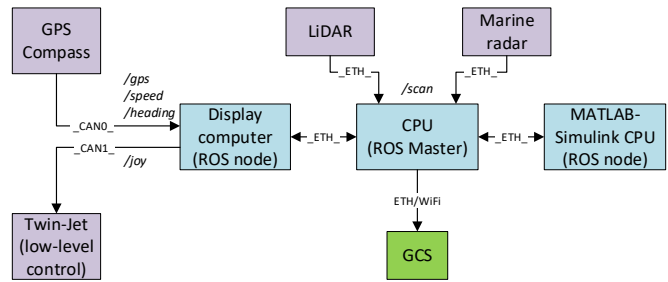

Fig. 3 Main USV GNC architecture

The model-based path planning and obstacle avoidance architecture for the twin jet USV (Villa, et al., 2019) uses the line-of-sight (LOS) path algorithm for straight-line as a GNC method. Despite the simple path following algorithm and architecture implemented in the USV, the error between the predefined and the field-test path is less than one meter. This difference is due to environmental elements (wind or wave drift forces), whose compensation is not considered in this study.

*2.6.3. Challenges for the future*

The described architecture used in GNC algorithms enables the implementation of different approaches, which can be developed in a simulation software (MATLAB/Simulink). To date, it has been tested with a simple path following algorithm (LOS path algorithm for a straight line). More complex algorithms will be implemented in the next phase (path planning with obstacle avoidance). This will make the architecture more complex, making it possible to have issues in the GNC implementation due to slow time-response of the system. The next step after testing the model-based architecture in the USV is to implement the same architecture in the rest of the offshore subsystems.

*2.7. Machine learning*

*2.7.1. Research goals*

The first step to ensure a safe cruise is to recognize the water surface (separate what is water and what is something else). The visual system signals to the USV control system that it is safe to steer the vessel forward if it detects open water alongside the planned route. Otherwise, the control system is advised to slow down, stop, or detour the USV to avoid a possible collision. The next step is to detect the objects in the water. Here, object recognition methods can be used to differentiate between diverse objects, such as rocks, boats (different in types and sizes), maritime signs, and other floating objects. Furthermore, it is aimed to predict the behaviour and paths of the moving objects. Then, the control system can apply different strategies when encountering various objects.

*2.7.2. Methods and results*

Our water segmentation approach is based on the KittiSeg road segmentation algorithm (Teichman, et al., 2016), which performs segmentation of the roads using a Fully Convolutional Network (FCN) deep learning architecture (Long, et al., 2015). The original model won the first place on the Kitti Road Detection Benchmark (Fritsch, et al., 2013) at the submission time. An important advantage of the KittiSeg algorithms is that the model can be trained using a small number of training samples. The original model used for road segmentation was trained with only 250 densely labelled images and achieved the state-of-the-art MaxF1 score of over 96%. With an inference time of 95ms per image using a GPU system, the model is suitable for real-time applications.

From the videos collected during the test runs, 600 images that have decent coverage of different sceneries and different types of objects visible in the image were selected. We annotated the selected images using the VGG Image Annotation (VIA) tool by marking the water area using polygons (Dutta, et al., 2016). To train our water segmentation model, we applied a transfer learning technique: we first initialized the weights of the FCN network using the weights of a pre-trained KittiSeg model and then we fine-tuned the model using the 600 manually labelled images. This transfer learning approach is extensively used in deep learning training to improve the accuracy of the model and reduce the training time. The resulting model was capable of recognizing the water surface successfully (see Fig. 4).

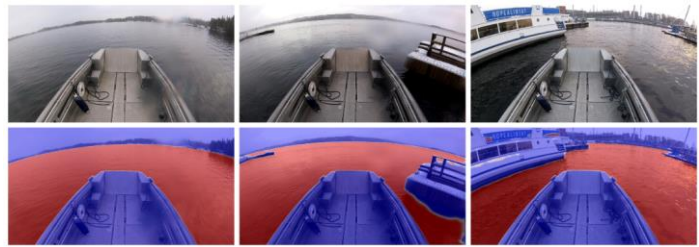

Fig. 4 Results of the water segmentation. Top row: the input image. Bottom row: the segmented image.

*2.7.3. Challenges for the future*

The next step after successfully performing the water segmentation is to detect and recognize objects in the water. The image classification and object detection are well-known computer vision applications. Systems that are trained using ImageNet dataset (Russakovsky, et al., 2015) and COCO dataset (Lin, et al., 2014) have led to promising classification and detection results. Our initial approach will be to fine-tune an existing model to perform object recognition on the water surface. Collecting additional data requires the use of the research vessel in real-world conditions, while labelling them remains a tedious and time demanding process. After obtaining sufficient labelled data, a model can be trained by using similar transfer learning techniques as before. More advanced techniques can combine the information received from these various USV sensors by applying multi-view and multi-modal learning for water segmentation and object recognition (Muslea, et al., 2002), (Ngian, et al., 2011). Multi-view learning permits to successfully recognize an image from different angles. The latter can be obtained by having one view from the USV and other views from the UAVs flying above it. In the concluding phase, we will aim to predict movements and behaviour of the observed objects, since predicting the paths of other vessels nearby and understanding how the maritime signs will affect their paths is essential for a smooth operation of the target USV.

*2.8. Communications*

*2.8.1. Research goals*

On the communication side, the project's main goal is to enable inexpensive low-latency, high-speed, and long-range bidirectional connection between all the autonomous subsystems. It is important to understand that today's communication systems typically target a subset of those characteristics. For example, satellite communications have broad coverage but high latency and unbalanced throughput performance for uplink and downlink (Peyravi, 1999). For our autonomous system, all of these metrics are important, since it is required to not only monitor the situation during the USV/UAV/AUV movement, but also to be able to interrupt the mission or permit for manual control of the subsystems in case of emergency.

*2.8.2. Methods and results*

Fig. 5 displays the layout of the general communication topology. It includes all offshore vehicles (USV, AUV, and UAV) as well as the GCS, which is mainly used for mission planning and control tasks. During the first phase of the project, we concentrated on a communication link between the GCS and the USV.

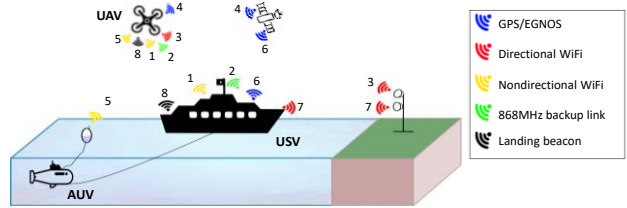

Fig. 5 General layout of aColor communication.

In order to achieve long-range and high-throughput communication, we decided to modify commercially available solutions, which are targeted to enable high-speed connections over a long distance assuming static positions of the transmitting and receiving sides. In particular, we chose IEEE 802.11ac (WiFi) directional interface DynaDish 5 produced by Mikrotik as our main communication device installed on both GCS and USV sides, and 3GPP LTE interface serving as a backup connection. The applied modifications enable mechanical rotation of the antenna on both GCS and USV sides controlled by a single-board computer (BeagleBone). On the software side, BeagleBone is running an algorithm that receives the global coordinates and orientation of the USV from the GPS compass and rotates the USV antenna towards the location of the GCS. The position of the ground station is assumed to be fixed and known to the USV. The GCS is running a similar algorithm where the position of the USV is received via an established communication link.

The main challenges encountered during our measurement campaign included the uncertainty of the GPS compass measurements and imperfections of the installed mechanical system (motors). In Fig. 6(a), the travel route and the appropriate RSSI measurements for WiFi and LTE interfaces are shown. The position points on the x-axis of the right subfigure of Fig. 6 correspond to the route steps shown in the left one. As it can be observed in Fig. 6(b), the system preserves acceptable WiFi RSSI using the described rotation mechanism, until the LOS path is blocked by an island. On the LTE link, the RSSI is stable all the time due to better coverage. However, LTE has much lower achievable throughput in comparison with 802.11ac (measured 400Mbps on 802.11ac vs. 100Mbps on LTE), and potentially exhibits higher latency because the data will first travel through the core network before being received on the GCS side.

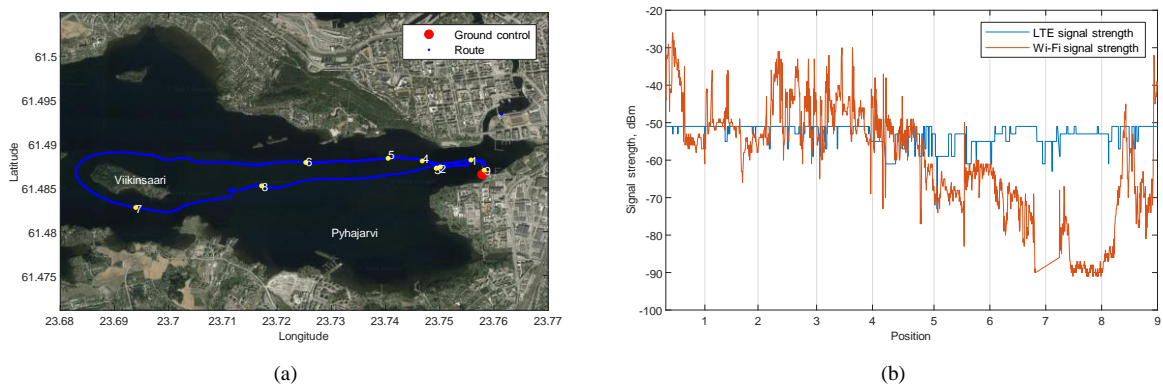

Fig. 6 Demo route (a) and corresponding RSSI measurements (b)

*2.8.3. Challenges for the future*

Further communications related tasks include testing of the UAV as a communication unit, which can work as a relay link between the GCS and the USV or transmit the collected data to one of those subsystems for subsequent processing. The first case becomes relevant when the USV-GCS communication link becomes occluded by an obstacle, e.g. an island. In that scenario, the UAV can rise to a certain altitude and relay the data, thus preserving high throughput and low delay on the appropriate links. However, this requires two high-speed interfaces utilizing the same communication technology and frequencies to be installed on the UAV, where at least one of the interfaces should have a similar range as the GCS or the USV directional transceivers. In our implementation, we employ Mikrotik SXT 5 ac, which has slightly lower range than the equipment installed on the GCS and the USV but it also has much smaller physical dimensions and weight, which allows us to install it on the UAV. Finally, the second interface used to connect UAV and USV is omnidirectional and does not require any mechanical support.

## 3. Other Autonomous Offshore Vehicles for cooperative systems

In order to enhance the USV reliability and robustness to system failures, improve mission performance, and optimize our strategies for a broader coverage of surveillance, communication, and measurement applications, it is essential to design and deploy a system that goes beyond a standalone USV.

*3.1. AUV and UAV architecture*

The aColor AUV is able to improve the capabilities of the complete offshore system by including the underwater areas into collision avoidance and object recognition considerations (using a camera and a mechanical sonar). This vehicle type can be used for underwater inspection to complement the 3D point-cloud created by the USV and UAV. In addition, ROS is used as a framework thus sharing the same interface as the rest of the vehicles of the offshore system. Fig. 7(a) outlines the AUV architecture.

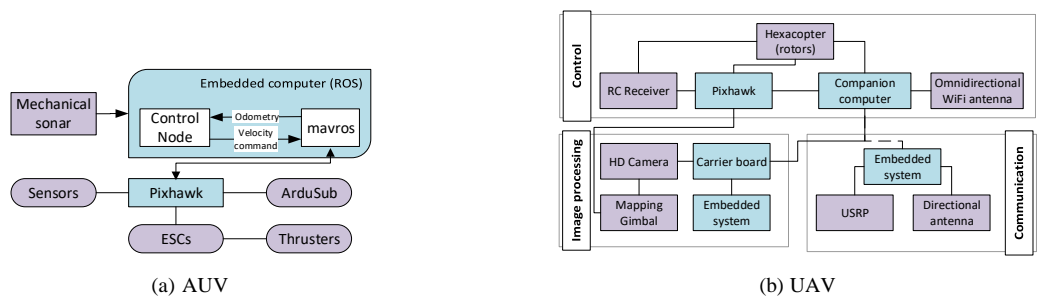

(a) AUV  (b) UAV

Fig. 7 Other offshore vehicles architecture

The aColor UAV is a multilayer system that enables different applications depending on the research topic being addressed. One of the main matters to be researched is to explore the application of moving relay stations (moving networks) for resilient connectivity by employing the UAV as a distributed antenna (i.e., relay link between the USV and the AUV). This way, the connection between GCS and USV can be improved by adding an extra communication link. Furthermore, the UAV can be utilized to provide improved obstacle avoidance capabilities. Hence, the UAV can be used mainly for search and rescue missions or as a moving relay station. Fig. 7(b) shows the UAV architecture with essential instrumentation that can be employed depending on the aerial vehicle.

*3.1.1. Research goals*

An offshore scenario differs remarkably from the regular inland UAV usage conditions because the USV is not stable in its position but can continuously move in the direction of all three axes. In addition, there is usually more wind to cope with than in normal landing onto the ground. The main research question set hereby is then: "How to move an UAV autonomously toward its charging position while landing it on the deck of an USV?"

*3.1.2. Methods and results*

The design of an UAV charging procedure includes the development of an electromechanical solution for an UAV to crib to a surface, to move to an accurate position required, to connect the battery to the charger, to start charging, and finally to return to fly. As a development method, Set-Based Design (SBD) (Sobek, et al., 1999) was applied to propose different solutions for each of the charging phases. With the SBD, it is possible to switch between the alternative solutions almost at the end of the product development process, if any loop-backs are needed later on (Camarda, et al., 2019). The two most challenging tasks in the automated charging process identified so far are the cribbing of the landing UAV to the landing platform and the moving of the UAV on the landing platform surface to the correct position. For both tasks, different solutions were sought and evaluated for implementation purposes.

For the cribbing, the following two solutions were considered:
1A. Use of reversing propellers to induce down force: While reversing, due to aerofoil, the fixed-pitch propellers produce less force than rotating normally. However, the force is enough for cribbing purposes.
1B. Use of electro permanent magnet gripper: The magnet type used was OpenGrab EPM v3, which has typical holding force of 300 N and weighs 36 g (NicaDrone, 2019). The magnet can be used to attach the UAV on the deck immediately after landing and before the mechanical arm comes to hold it.

To centre the UAV on the landing platform surface at the USV includes:
2A. Use of linear electrical actuators: This method was found easy to implement and robust against different weather conditions.
2B. Use of gravity to centre the drone on the landing platform: This method is simpler and more robust because it avoids moving components in the landing platform.

As our initial solutions in this study, the most robust ones for harsh conditions were selected from the aforementioned competing options. With the selections done, the entire automated landing and charging procedure is as follows: the UAV lands and applies the magnet to crib the landing gear when it reaches the bottom of the platform. Gravity and the shape of the landing gear and the landing platform are used as basis to move the UAV to the correct charging position. The inductive sensors on the UAV's landing gear indicate the touch down and switch off the UAV motors, sending this information to the landing platform. The charging connector moves by its own actuator and the mechanical shape, while a spring-based structure enables adequate electrical connection. After the charging is complete, the connector moves out and a message is sent to the UAV. Finally, the UAV starts its motors, demagnetises the magnet, and starts to fly.

*3.1.3. Challenges for the future*

Due to a lack of commercial solutions and literature, the SBD method was applied to outline several competing solutions, of which the most robust for harsh conditions was selected to implement and test. The electro permanent magnet gripper and linear electrical actuators with flexible charging connector are proposed as the preferred initial solution to be used. However, this approach has not yet been deployed in a real scenario. In our future studies, the use of different criteria or requirements and field test results comparing the competing solutions should help understand the steps towards the final solution.

*3.2. Communication between offshore vehicles (landing system)*

*3.2.1. Research goals*

An automated UAV terminal guidance and recovery system is needed on the USV to ensure reliable landing of the UAV in all weather conditions. In this study, a combination of visual and Radio Frequency (RF) methods is used to facilitate the landing process. The landing system is designed to operate in adverse weather conditions (visibility below 30 meters, rain or snow, high wind). To this end, the landing system follows staged approach, incrementally switching to more precise (and, consequently, shorter range) systems.

*3.2.2. Methods*

The necessary equipment in the USV features the landing platform itself that performs mechanical capture of the UAV and WiFi RF sensing antennas used to measure the signal strength on the link towards the UAV. All antennas are connected to the same receiver, which permits precise measurements of the signal power and thus accurate localization. In addition to that, coloured LED lamps (red, green, and blue) are setup to shine vertically up for the UAV to detect once it is hovering on top of the USV. For the UAV, the specific equipment for the landing procedure consists of an omnidirectional WiFi antenna for the signal strength measurements, a down-facing low-resolution camera for the detection of landing lights, an ultrasonic altimeter to determine the exact distance to the landing pad during the final approach, and the electro permanent magnet gripper to lock the UAV on the landing platform once a touchdown is achieved.

## 4. Conclusions

In conclusion, the research areas of mechatronics, machine learning, and communications need to be intertwined to produce a successful autonomous offshore multicomponent robotic system. From the mechatronics point of view, model-based architecture allows for a rapid-prototyping procedure for the GNC algorithms in any of the offshore subsystems. From the machine learning side, it is possible to include the obstacle avoidance capabilities based on a water segmentation algorithm, which permits the detection and recognition of objects in the water. From the communication side, inexpensive low-latency, high-speed, and long-range bidirectional connections between autonomous subsystems offer the possibility to monitor and control the subsystems in case of emergency. Regarding cooperation between different offshore subsystems, the main challenges and approaches for providing a robust mechanical solution to move and restrain a landing UAV to a predefined position on a surface for battery charging are also presented. Finally, the methods for accurate landing of the UAV in the USV platform are presented using a combination of visible and RF cues, which ensures operability in all weather conditions.

**Acknowledgment**

This paper is based on the aColor project funded by Technology Industries of Finland Centennial and Jane & Aatos Erkko Foundations under Future Makers Funding Program 2017. The authors gratefully acknowledge the contributions of the company Alamarin-Jet Oy.

**References**

Camarda, D. J., Scotti, S., Kunttu, I. & Perttula, A., 2019. Rapid product development methods in practice – case studies from industrial production and technology development. Proceedings of ISPIM Connects.

Dutta, A., Gupta, A. & Zissermann, A., 2016. VGG image annotator (VIA).. [Online] Available at: http://www. robots. ox. ac. uk/~vgg/software/via [Accessed 12 October 2019].

Eriksson, S., 2016. Real-time kinematic positioning of UAS - possibilities and restrictions, Gothenburg, Sweden: Master's thesis. Chalmers University of Technology.

Fossen, T. I., Breivik, M. & Skjetne, R., 2003. Line-of-sight path following of underactuated marine craft.. IFAC Proceedings Volumes 36, no. 21, pp. 211-216.

Fritsch, J., Kuehnl, T. & Geiger, A., 2013. A new performance measure and evaluation benchmark for road detection algorithms.. International Conference on Intelligent Transportation Systems (ITSC).

Lin, T.-Y.et al., 2014. Microsoft coco: Common objects in context.. European conference on computer vision, pp. 740-755.

Long, J., Shelhamer, E. & Darrell, T., 2015. Fully convolutional networks for semantic segmentation.. Proceedings of the IEEE Conference on Computer Vision and Pattern Recognition, p. 3431–3440.

Moe, S. & Pettersen, K. Y., 2016. Set-based Line-of-Sight (LOS) path following with collision avoidance for underactuated unmanned surface vessel.. 24th Mediterranean Conference on Control and Automation (MED), IEEE, pp. 402-409.

Muslea, I., Minton, S. & Knoblock, C. A., 2002. Active+ semi-supervised learning = robust multi-view learning.. Proceedings of the International Conference on Machine Learning..

Ngian, J. et al., 2011. Multimodal deep learning.. Proceedings of the 28th international conference on machine learning (ICML-11), pp. 689-696.

NicaDrone, 2019. OpenGrab EPM V3 R5C, datasheet. [Online] Available at: https://nicadrone.com/ [Accessed 12 October 2019].

Peyravi, H., 1999. Medium access control protocols performance in satellite communications. IEEE Communications Magazine 37.3, pp. 62-71.

Quigley, M. et al., 2019. ROS: an open-source Robot Operating System.. Kobe, Japan, ICRA workshop on open source software..

Russakovsky, O. et al., 2015. ImageNet Large Scale Visual Recognition Challenge. IJVC.

Sobek, D., Ward, A. C. & Liker, J. K., 1999. Toyota's Principles of Set-Based Concurrent Engineering. Sloan Management Review, Volume 40 (2), p. 67–83.

solidThinking Inspire, 2018. Topology Optimization. [Online] Available at: https://solidthinking.com/help/Inspire/2018/win/en_us/topology_optimization.htm [Accessed 12 October 2019].

TechFinland100, 2018. Autonomous and Collaborative Offshore Robotics (aCOLOR),. [Online] Available at: https://techfinland100.fi/mita-rahoitamme/tutkimus/tulevaisuuden-tekijat/autonomous-and-collaborative-offshore-robotics-acolor/ [Accessed 12 October 2019].

Teichman, M. et al., 2016. MultiNet: Real-time Joint Semantic Reasoning for Autonomous Driving. CoRR.

The Mathworks, Inc., 2018. Generate a Standalone ROS Node from Simulink. [Online] Available at: https://se.mathworks.com/help/robotics/examples/generate-a-standalone-ros-node-in-simulink.html [Accessed 12 October 2019].

Turkkila, T. & Ukonaho, M., 2016. Study About Boundaries and Constraints in Topology Optimization. NAFEMS Event: Exploring the Design Freedom of Additive Manufacturing through Simulation.

Villa, J., Aaltonen, J. & Koskinen, K. T., 2019. Model-Based path planning and obstacle avoidance Architecture for a Twin Jet Unmanned Surface Vessel. Third IEEE International Conference on Robotic Computing (IRC), pp. 427-428.